\documentclass[11pt]{article}
\usepackage{amsmath}
\usepackage{amssymb}
\usepackage{color}
\usepackage{tikz}
\usepackage{endnotes}
\usepackage{biblatex}
\usepackage{graphicx}
\usepackage{colortbl}
\begin{document}
\title{On the relationship between the Collatz conjecture and Mersenne prime numbers}
\author{Jonas Kaiser\footnote{jonaskaiser79@gmx.de}}
\date{\today}
\maketitle
\begin{abstract}
The purpose of this study is to show how to get a necessary criterion for prime numbers with the help of special matrices. My special interest lies in the empirical research of these matrices and their patterns, structures and symmetries. The matrices in turn depend on an expansion of the Collatz algorithm 3n+1. 
\end{abstract}
\newpage
\tableofcontents
\newpage
\section{Introduction}
The basis of this study is an expansion of the Collatz conjecture. "The Collatz conjecture can be summarized as follows. Take any positive integer n. If n is even, divide it by 2 to get n / 2. If n is odd, multiply it by 3 and add 1 to obtain 3n + 1. Repeat the process indefinitely. The conjecture is that no matter what number you start with, you will always eventually reach 1." \cite{1}\\\\ 
For each expansion it is possible to assign a specific structure. The structures in turn can be identified by specific matrices, and these matrices have a relation to prime numbers (Collatz matrix conjecture). Furthermore, it can be shown that there is a special relation between matrices and Mersenne prime numbers (Kaiser's conjecture).\\\\
First of all in Chapter 2, i will expand the Collatzalgorithm and define a general form of a matrix based on these Collatzalgorithms. This so-called Collatz tree matrix can be identified by other special matrices (The Collatz matrices), i will define in Chapter 3. In Chapter 4 it can be shown how the Collatz matrices deliver a necessary criterion for prime numbers (Collatz matrix conjecture), and Chapter 5 contains a conjecture about Collatz matrices and Mersenne prime numbers (Kaiser's conjecture).\\\\
The definitions of Chapter 2 (2.1, 2.2, 2.3 and 2.4) are already well known to the "Collatz community". For further Informaton see: \cite{2,3}. All other definitions and conjectures were found by my own, if not explicitly advertised.
\section{The Collatz tree matrix} 
	\subsection{Expansion of the Collatz algorithm}
First i will expand the Collatzalgorithm.
\\\\
\begin{equation}
   f(n) =
   \begin{cases}
     n/2,  \ if \ n \ is \ even  \\        
     3n+1,  \ if \ n \ is \ odd 
   \end{cases}
\end{equation}
\\\\
with respect to a
\\\\
\begin{equation}
   f_a(n) =
   \begin{cases}
     n/2,  \ if \ n \ is \ even  \\        
     an+1,  \ if \ n \ is \ odd 
   \end{cases}
\end{equation}
\\\\
$$D=\{a|a=2n+1, \ n\in\mathbb{N}_0\}$$
	\subsection{The Collatz tree matrix}
\textbf{Definition 2.1 The Collatz tree matrix}\\
Each Collatzalgorithm $f_a(n)$ can be shown as a Collatz tree matrix. The Collatz tree matrix consists of three different types of columns. The first column is called the knot column. In this column all results of the specific odd operation are listed. The second column is called the odd column. In this column all natural odd numbers are listed. The third and each additional column is called the even column because the numbers of these columns are results of the even operation. The Collatz tree matrices only differ from their first column.\\\\
The general form of the Collatz tree matrix of $f_a(n)$ looks like that:\\\\
\begin{center} Collatz tree matrix of algorithm $f_a(n)$:
\end{center}
$$
\begin{pmatrix}
b_1a+1 &b_1 & 2b_1 & 2^2b_1 & 2^3b_1 & 2^4b_1 &\ldots &2^nb_1 \\
b_2a+1 & b_2 & 2b_2 & 2^2b_2 &2^3b_2&2^4b_2& \ldots &2^nb_2 \\
b_3a+1 & b_3 & 2b_3 & 2^2b_3&2^3b_3 &2^4b_3 & \ldots &2^nb_3\\
b_4a+1 & b_4 & 2b_3 & 2^2b_4&2^3b_4 &2^4b_4 & \ldots &2^nb_4\\
\vdots & \vdots & \vdots   & \vdots& \vdots  & \vdots  &\ddots & \vdots \\
b_na+1 & b_n & 2b_n & 2^2b_n&2^3b_n &2^4b_n & \ldots &2^nb_n \\
\end{pmatrix}
$$
\\
$$D=\{b|b=2n+1,n\in\mathbb{N}_0\}$$
\\\\
Let us have a look at the first three Collatz tree matrices:\\\\
\begin{center} Collatz tree matrix of $f_3(n)$:
\end{center}
$$
\begin{pmatrix}
4 & 1 & 2 & 4 & 8 & \ldots & \infty\\
10 & 3 & 6 & 12 & 24&\ldots & \infty\\
16 & 5 & 10 & 20 & 40 & \ldots& \infty\\
22& 7 & 14 & 28 & 56 & \ldots& \infty\\
28&9 & 18 & 36 & 72 & \ldots & \infty\\
34& 11 & 22 & 44 & 88 & \ldots & \infty\\
 \vdots & \vdots & \vdots & \vdots & \vdots & \ddots & \ \\
 \infty & \infty & \infty & \infty & \infty & \ & \infty\\
\end{pmatrix}
$$
\\\\
\begin{center} Collatz tree matrix of $f_5(n)$:
\end{center}
$$
\begin{pmatrix}
6 & 1 & 2 & 4 & 8 & \ldots & \infty\\
16 & 3 & 6 & 12 & 24 & \ldots & \infty\\
26 & 5 & 10 &20 & 40 & \ldots & \infty\\
36 & 7 & 14 & 28 & 56 & \ldots & \infty\\
46 & 9 & 18 & 36 & 72 &\ldots & \infty\\
56 & 11 & 22 & 44 & 88 & \ldots & \infty\\
 \vdots  & \vdots  & \vdots  & \vdots  &\vdots  & \ddots & \  \\
\infty & \infty & \infty & \infty & \infty & \ & \infty  \\
\end{pmatrix}
$$
\\\\
\begin{center} Collatz tree matrix of $f_7(n)$:
\end{center}
$$
\begin{pmatrix}
8 & 1 & 2 & 4 & 8 & \ldots & \infty\\
22 & 3 & 6 & 12 & 24 & \ldots & \infty\\
36 & 5& 10 & 20 & 40 & \ldots & \infty\\
50 & 7 & 14 & 28 & 56 & \ldots & \infty\\
64 & 9 & 18 & 36 & 72 &\ldots & \infty\\
78 & 11 & 22 & 44 & 88 & \ldots & \infty\\
 \vdots  & \vdots  & \vdots  & \vdots  & \vdots  & \ddots & \  \\
\infty & \infty & \infty & \infty & \infty & \ & \infty  \\
\end{pmatrix}
$$
\\\\
	\subsection{Knot numbers}
\textbf{Definition 2.2 Knot numbers}\\
Each Collatz tree matrix has knot numbers. They are even and result from the specific odd operation. Each knot number appears two times. One time they are listed in the first knot column of the Collatz tree matrix and the other time they appear right of the odd column. If a knot number appears two times in the same row, then there is a cycle. \cite{2,3}
\\\\
As can be seen each Collatz tree matrix has his own specific knot number pattern (red numbers):
\\\\
\begin{center} Collatz tree matrix of $f_3(n)$:
\end{center}
$$
\begin{pmatrix}
\textcolor{red}{4} & 1 & 2 & \textcolor{red}{4} &8&\ldots& \infty\\
\textcolor{red}{10} & 3 & 6 & 12 & 24& \ldots &\infty\\
\textcolor{red}{16} & 5 & \textcolor{red}{10} & 20 & \textcolor{red}{40} & \ldots &\infty\\
\textcolor{red}{22} & 7 & 14 & \textcolor{red}{28} & 56 &\ldots &\infty\\
\textcolor{red}{28} & 9 & 18 & 36 & 72 &\ldots &\infty\\
\textcolor{red}{34} & 11 & \textcolor{red}{22} & 44 & \textcolor{red}{88} & \ldots&\infty\\
 \vdots & \vdots & \vdots & \vdots &\vdots & \ddots &  \ \\
 \infty & \infty & \infty & \infty & \infty & \ & \infty \\
\end{pmatrix}
$$
\\\\
\begin{center} Collatz tree matrix of $f_5(n)$:
\end{center}
$$
\begin{pmatrix}
\textcolor{red}{6} & 1 & 2 & 4 & 8 &  \textcolor{red}{16} & \ldots &\infty\\
\textcolor{red}{16} & 3 & \textcolor{red}{6} & 12 & 24 & 48 & \ldots  &\infty\\
\textcolor{red}{26} &5 & 10 & 20 & 40 & 80 & \ldots & \infty\\
\textcolor{red}{36} & 7 & 14 & 28 & \textcolor{red}{56} & 112 &\ldots  &\infty\\
\textcolor{red}{46} & 9 & 18 & \textcolor{red}{36} & 72 & 144 &\ldots & \infty\\
\textcolor{red}{56} & 11 & 22 & 44 & 88 & \textcolor{red}{176} & \ldots  &\infty\\
\textcolor{red}{66} & 13 & \textcolor{red}{26} & 52 & 104 & 208 & \ldots & \infty\\
\textcolor{red}{76} & 15 & 30 &60 &120 & 240 & \ldots  &\infty\\
\vdots & \vdots & \vdots & \vdots &\vdots &\vdots& \ddots & \ \\
\infty & \infty & \infty & \infty & \infty & \infty& \ &\infty  \\
\end{pmatrix}
$$
\\\\
\begin{center} Collatz tree matrix of $f_7(n)$:
\end{center}
$$
\begin{pmatrix}
\textcolor{red}{8} & 1 & 2 & 4 & \textcolor{red}{8} & \ldots  &\infty\\
\textcolor{red}{22} & 3 & 6 & 12 & 24 & \ldots  &\infty\\
\textcolor{red}{36} &5 & 10 & 20 & 40 & \ldots &\infty\\
\textcolor{red}{50} & 7 & 14 & 28 & 56 & \ldots  &\infty\\
\textcolor{red}{64} & 9 & 18 & \textcolor{red}{36} & 72 & \ldots  &\infty\\
\textcolor{red}{78} & 11 & \textcolor{red}{22} & 44 & 88 & \ldots &\infty\\
\textcolor{red}{91} & 13 & 26 & 52 & 104 & \ldots & \infty\\
\vdots & \vdots & \vdots & \vdots &\vdots &\ddots &\ \\
\infty & \infty & \infty & \infty & \infty & \ &\infty \\
\end{pmatrix}
$$
	\subsection{Unbranched rows}
\textbf{Definition 2.3 Unbranched rows}\\
Each Collatz tree matrix has unbranched rows. Such a row has only one knot number. \cite{2,3}
\\\\
They look like that (blue rows):\\\\
\begin{center} Collatz tree matrix of $f_3(n)$:
\end{center}
$$
\begin{pmatrix}
\textcolor{red}{4} & 1 & 2 & \textcolor{red}{4} & 8 & \ldots &\infty\\
\textcolor{red}{10} & \textcolor{blue}{3} & \textcolor{blue}{6} & \textcolor{blue}{12} & \textcolor{blue}{24}& \ldots &\infty\\
\textcolor{red}{16} & 5 & \textcolor{red}{10} & 20 & \textcolor{red}{40} & \ldots&\infty\\
\textcolor{red}{22} & 7 & 14 & \textcolor{red}{28} & 56 & \ldots &\infty\\
\textcolor{red}{28} & \textcolor{blue}{9} & \textcolor{blue}{18} & \textcolor{blue}{36} & \textcolor{blue}{72} & \ldots &\infty\\
\textcolor{red}{34} & 11 & \textcolor{red}{22} & 44 & \textcolor{red}{88} &\ldots&\infty\\
\vdots & \vdots & \vdots & \vdots &\vdots &\ddots &\ \\
 \infty & \infty & \infty & \infty & \infty & \ & \infty\\
\end{pmatrix}
$$
\\\\
\begin{center} Collatz tree matrix of $f_5(n)$:
\end{center}
$$
\begin{pmatrix}
\textcolor{red}{6} & 1 & 2 & 4 & 8 &  \textcolor{red}{16} & \ldots&\infty\\
\textcolor{red}{16} & 3 & \textcolor{red}{6} & 12 & 24 & 48 & \ldots&\infty\\
\textcolor{red}{26} & \textcolor{blue}{5} & \textcolor{blue}{10} & \textcolor{blue}{20} & \textcolor{blue}{40} & \textcolor{blue}{80} &\ldots &\infty\\
\textcolor{red}{36} & 7 & 14 & 28 & \textcolor{red}{56} & 112 &\ldots& \infty\\
\textcolor{red}{46} & 9 & 18 & \textcolor{red}{36} & 72 & 144 & \ldots&\infty\\
\textcolor{red}{56} & 11 & 22 & 44 & 88 & \textcolor{red}{176} &\ldots&\infty\\
\textcolor{red}{66} & 13 & \textcolor{red}{26} & 52 & 104 & 208 & \ldots&\infty\\
\textcolor{red}{76} & \textcolor{blue}{15} & \textcolor{blue}{30} & \textcolor{blue}{60} & \textcolor{blue}{120} & \textcolor{blue}{240} & \ldots &\infty\\
\vdots & \vdots & \vdots & \vdots &\vdots & \vdots &\ddots &\ \\
\infty & \infty & \infty & \infty & \infty & \infty & \ & \infty \\
\end{pmatrix}
$$
\\\\
\begin{center} Collatz tree matrix of $f_7(n)$:
\end{center}
$$
\begin{pmatrix}
\textcolor{red}{8} & 1 & 2 & 4 & \textcolor{red}{8} & \ldots&\infty\\
\textcolor{red}{22} & \textcolor{blue}{3} & \textcolor{blue}{6} & \textcolor{blue}{12} & \textcolor{blue}{24} & \ldots&\infty\\
\textcolor{red}{36} & \textcolor{blue}{5} & \textcolor{blue}{10} & \textcolor{blue}{20} & \textcolor{blue}{40} & \ldots &\infty\\
\textcolor{red}{50} & \textcolor{blue}{7} & \textcolor{blue}{14} & \textcolor{blue}{28} & \textcolor{blue}{56} & \ldots&\infty\\
\textcolor{red}{64} & 9 & 18 & \textcolor{red}{36} & 72 & \ldots &\infty\\
\textcolor{red}{78} & 11 & \textcolor{red}{22} & 44 & 88 & \ldots&\infty\\
\textcolor{red}{91} & \textcolor{blue}{13} & \textcolor{blue}{26} & \textcolor{blue}{52} & \textcolor{blue}{104} & \ldots&\infty\\
\vdots & \vdots & \vdots & \vdots &\vdots &\ddots &\ \\
\infty & \infty & \infty & \infty & \infty & \ & \infty \\
\end{pmatrix}
$$
	\subsection{Perfect knot numbers}
\textbf{Definition 2.4 Perfect knot numbers}\\
Perfect knot numbers are rooted in unbranched rows. \cite{2,3}
\\\\
Let us have a look at the perfect knot numbers (frameboxed numbers) of the Collatz tree matrix of $f_3(n)$:\\\\
\\
\begin{center}Collatz tree matrix of $f_3(n)$
\end{center}
$$
\begin{pmatrix}
\textcolor{red}{4}&1&2 & \textcolor{red}{4}& 8 & \textcolor{red}{16} & 32 & \framebox[8mm]{\textcolor{red}{64}}  & \ldots &\infty\\
\framebox[8mm]{\textcolor{red}{10}}↨&3&\textcolor{blue}{6} & \textcolor{blue}{12} & \textcolor{blue}{24} & \textcolor{blue}{48} & \textcolor{blue}{96} & \textcolor{blue}{192} &  \ldots& \infty\\
\textcolor{red}{16}&5&\framebox[8mm]{\textcolor{red}{10}} & 20 & \textcolor{red}{40} & 80 & \textcolor{red}{160} & 320  & \ldots &\infty\\
\textcolor{red}{22}&7&14 & \framebox[8mm]{\textcolor{red}{28}} & 56 & \textcolor{red}{112} & 224 & \textcolor{red}{448} &  \ldots&\infty \\
\framebox[8mm]{\textcolor{red}{28}}&9&\textcolor{blue}{18} & \textcolor{blue}{36} & \textcolor{blue}{72} & \textcolor{blue}{144} & \textcolor{blue}{288} & \textcolor{blue}{576}  & \ldots&\infty\\
\textcolor{red}{34}&11&\textcolor{red}{22} & 44 & \textcolor{red}{88} & 176 & \framebox[8mm]{\textcolor{red}{352}} & 704  & \ldots&\infty\\
\textcolor{red}{40}&13&26 & \textcolor{red}{52} & 104 & \framebox[8mm]{\textcolor{red}{208}} & 416 & \textcolor{red}{832}  &  \ldots&\infty\\
\framebox[8mm]{\textcolor{red}{46}}&15&\textcolor{blue}{30} & \textcolor{blue}{60} & \textcolor{blue}{120} & \textcolor{blue}{240} & \textcolor{blue}{480} & \textcolor{blue}{960} & \ldots& \infty \\
\textcolor{red}{52}&17&\textcolor{red}{34} & 68 & \framebox[8mm]{\textcolor{red}{136}} & 272 & \textcolor{red}{544} & 1088 &  \ldots & \infty \\
\vdots & \vdots & \vdots & \vdots &\vdots & \vdots & \vdots &\vdots &\ddots&\ \\
\infty & \infty & \infty & \infty & \infty & \infty & \infty & \infty &  \ & \infty \\
\end{pmatrix}
$$
\\\\
\section{The Collatz matrix}
	\subsection{The three Collatz matrices}
With these knot numbers, unbranched rows and perfect knot numbers each Collatz tree matrix can be identified by three Collatz matrices i will define in the next sections. They are called \textit{Standard Collatz matrix}, \textit{Little Collatz matrix} and \textit{Big Collatz matrix}. 
	\subsection{The standard Collatz matrix}
\textbf{Definition 3.1 The standard Collatz matrix}\\
Each Collatz tree matrix can be identified by a specific standard Collatzmatrix. The standard Collatz matrix is the area (of the Collatz tree matrix), whose structure (pattern of knot numbers and unbranched rows) repeats. $a_{11}$ of the standard Collatz matrix is always $2b_1$ of the general Collatz tree matrix. The number of rows of a standard Collatz matrix (of the Collatz tree matrix) of algorithm $f_a(n)$ is always $a$. The number of columns depends on the first knot number in the first row. In other words: $a_{1n}$ of the standard Collatz matrix is always the first knot number of that row.\\\\
The number of rows $m_C$ is in the following relationship to the number of columns $n_C$:\\\\
$$(2^{n_C}-1)\ mod \ m_C = 0$$
$$m_C>n_C$$
\\\\
The first three standard Collatz matrices are:
\\\\
\begin{center}Standard Collatz matrix 3 x 2:
\end{center}
$$
\begin{pmatrix}
 2 & \textcolor{red}{4} \\
 \textcolor{blue}{6} & \textcolor{blue}{12} \\
 \textcolor{red}{10} & 20 \\
\end{pmatrix}
$$
\\\\
\begin{center}Standard Collatz matrix 5 x 4:
\end{center}
$$
\begin{pmatrix}
 2 & 4& 8 & \textcolor{red}{16} \\
 \textcolor{red}{6} & 12 & 24 & 48\\
 \textcolor{blue}{10} & \textcolor{blue}{20} & \textcolor{blue}{40} & \textcolor{blue}{80} \\
 14 & 28 & \textcolor{red}{56} & 112 \\
 18 & \textcolor{red}{36} & 72 & 144 \\
\end{pmatrix}
$$
\\\\
\begin{center}Standard Collatz matrix 7 x 3:
\end{center}
$$
\begin{pmatrix}
 2 & 4& \textcolor{red}{8}  \\
 \textcolor{blue}{6} & \textcolor{blue}{12} & \textcolor{blue}{24} \\
 \textcolor{blue}{10} & \textcolor{blue}{20} & \textcolor{blue}{40} \\
 \textcolor{blue}{14} & \textcolor{blue}{28} & \textcolor{blue}{56} \\
 18 & \textcolor{red}{36} & 72 \\
 \textcolor{red}{22} & 44 & 88 \\
 \textcolor{blue}{26} & \textcolor{blue}{52} & \textcolor{blue}{104} \\
\end{pmatrix}
$$
\clearpage
To get a better view on the standard Collatz matrix let us have a look in context of the Collatz tree matrix (Figure 1):
\\\\
\begin{center}Standard Collatz matrix 3 x 2 (green box) as part of the
\end{center}
\begin{center}Collatz tree matrix of $f_3(n)$
\end{center}
\begin{figure}[htbp] 
  \centering
     \includegraphics[width=1\textwidth]{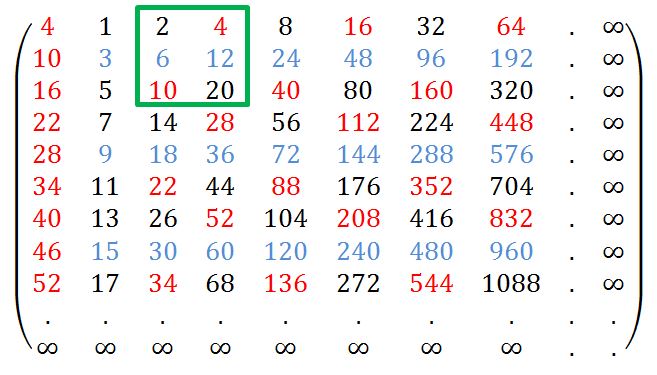}
  \caption{Collatz tree matrix of $f_3(n)$}
  \label{fig:Bild1}
\end{figure}
As can be seen the structure (position of knot numbers and unbranched rows) of the standard Collatz matrix (green box) repeats downwards and to the right.
\\\\
Important Notice: When i write of the Collatz matrix in the next sections i always mean the standard Collatz matrix.
\clearpage
	\subsection{Little Collatz matrix}
\textbf{Definition 3.2 The little Collatz matrix}\\
The little Collatz matrix is part of the Collatz matrix. It has the same value for $a_{11}$ like the Collatz matrix but it ends with the first (lowest) knot number.
\\\\
For example the Collatz tree matrix of $f_3(n)$ has the little Collatz matrix:
\\\\
\begin{center}
little collatz matrix 1 x 2: 
\end{center}
$$
\begin{pmatrix}
2 & \textcolor{red}{4} \\
\end{pmatrix}
$$
\\\\
The Collatz tree matrix of $f_5(n)$ has the little Collatz matrix:
\\\\
\begin{center}
Little collatz matrix 2 x 1: 
\end{center}
$$
\begin{pmatrix}
2 \\
\textcolor{red}{6} \\
\end{pmatrix}
$$
\\\\
The Collatz tree matrix of $f_{23}(n)$ has the little Collatz matrix:
\\\\
\begin{center}
Little Collatz matrix 2 x 3: 
\end{center}
$$
\begin{pmatrix}
2 & 4 & 8\\
6 & 12 & \textcolor{red}{24}\\
\end{pmatrix}
$$
\\\\
The little Collatz matrix is related to Collatz matrix like this :
\\\\
$$\frac{1+m_C+2^{n_L}}{2^{n_L+1}}=m_L$$
\\\\
$m_C=m$-value of the standard Collatz matrix\\\\
$m_L=m$-value of the little Collatz matrix\\\\
$n_L=n$-value of the little Collat zmatrix\\\\
	\subsection{Big Collatz matrix}
\textbf{Definition 3.3 The big Collatz matrix}\\
Each Collatz tree matrix can be identified by a specific big Collatz matrix. The big Collatz matrix is the area (of the Collatz tree matrix), whose structure (pattern of perfect knot numbers) repeats. $a_{11}$ of the big Collatz matrix is always $2b_1$ of the general Collatz tree matrix. The number of rows of a big Collatz matrix (of the Collatz tree matrix) of algorithm $f_a(n)$ is always $a^2$. 
\\\\
The big Collatz matrix of $f_3(n)$ with the perfect knot numbers (frameboxed numbers):
\\\\
\begin{center}Big Collatz matrix 9 x 6
\end{center}
$$
\begin{pmatrix}
2 & \textcolor{red}{4}& 8 & \textcolor{red}{16} & 32 & \framebox[8mm]{\textcolor{red}{64}} \\
\textcolor{blue}{6} & \textcolor{blue}{12} & \textcolor{blue}{24} & \textcolor{blue}{48} & \textcolor{blue}{96} & \textcolor{blue}{192}\\
\framebox[8mm]{\textcolor{red}{10}} & 20 & \textcolor{red}{40} & 80 & \textcolor{red}{160} & 320 \\
14 & \framebox[8mm]{\textcolor{red}{28}} & 56 & \textcolor{red}{112} & 224 & \textcolor{red}{448} \\
\textcolor{blue}{18} & \textcolor{blue}{36} & \textcolor{blue}{72} & \textcolor{blue}{144} & \textcolor{blue}{288} & \textcolor{blue}{576} \\
\textcolor{red}{22} & 44 & \textcolor{red}{88} & 176 & \framebox[8mm]{\textcolor{red}{352}} & 704 \\
26 & \textcolor{red}{52} & 104 & \framebox[8mm]{\textcolor{red}{208}} & 416 & \textcolor{red}{832} \\
\textcolor{blue}{30} & \textcolor{blue}{60} & \textcolor{blue}{120} & \textcolor{blue}{240} & \textcolor{blue}{480} & \textcolor{blue}{960} \\
\textcolor{red}{34} & 68 & \framebox[8mm]{\textcolor{red}{136}} & 272 & \textcolor{red}{544} & 1088 \\
\end{pmatrix}
$$
\\\\
The m- and n-values of the big Collatz matrix are given by that:\\\\
$$m_B \ x \ n_B=(m_C*m_C) x (m_C*n_C)$$\\\\
$m_B=m$-value of the big Collatz matrix\\\\
$n_B=n$-value of the big Collatz matrix\\\\
$m_C=m$-value of the Standard Collatz matrix\\\\
$n_C=n$-value of the Standard Collatz matrix\\\\
	\subsection{The three Collatz matrices as part of the Collatz tree matrix}
Here again pictures of the Collatz tree matrices of $f_3(n)$ (Figure 2),$f_5(n)$ (Figure 3) and $f_7(n)$ (Figure 4) with their three matrices (blue: Little Collatz matrix, green: Standard Collatz matrix and red: Big Collatz matrix):
\clearpage
\begin{figure}[htbp] 
  \centering
     \includegraphics[width=1\textwidth]{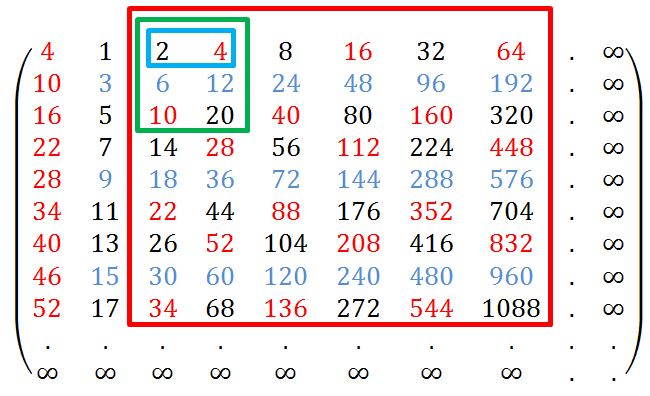}
  \caption{Collatz tree matrix of $f_3(n)$}
  \label{fig:Bild1}
\end{figure}
\begin{center}blue: Little Collatz matrix 1 x 2\\
green: Standard Collatz matrix 3 x 2\\
red: Big Collatz matrix 9 x 6
\end{center}
\clearpage
\begin{figure}[htbp] 
  \centering
     \includegraphics[width=1\textwidth]{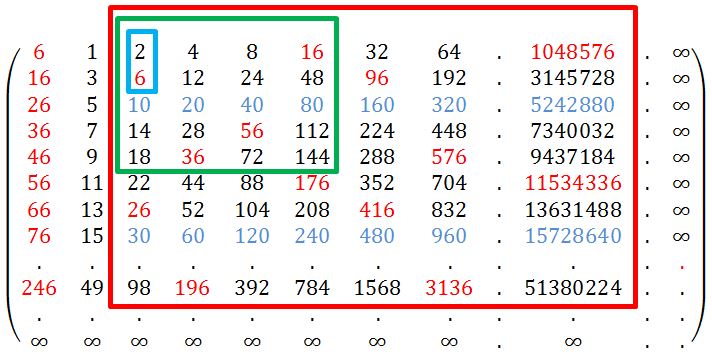}
  \caption{Collatz tree matrix of $f_5(n)$}
  \label{fig:Bild1}
\end{figure}
\begin{center}blue: Little Collatz matrix 2 x 1\\
green: Standard Collatz matrix 5 x 4\\
red: Big Collatz matrix 25 x 20
\end{center}
\clearpage
\begin{figure}[htbp] 
  \centering
     \includegraphics[width=1\textwidth]{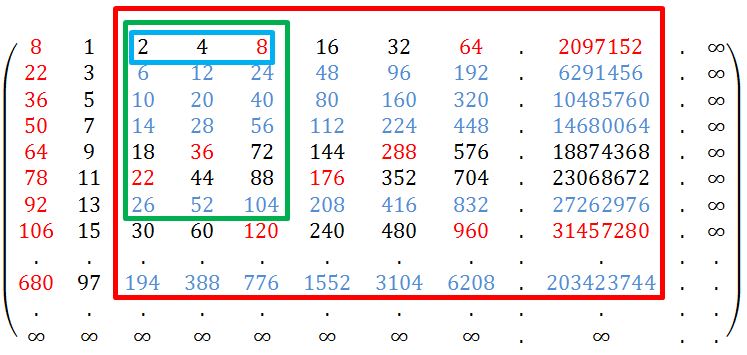}
  \caption{Collatz tree matrix of $f_7(n)$}
  \label{fig:Bild1}
\end{figure}
\begin{center}blue: Little Collatz matrix 1 x 3\\
green: Standard Collatz matrix 7 x 3\\
red: Big Collatz matrix 49 x 21
\end{center}
\clearpage
	\subsection{Symmetries and structures}
When comparing the standard Collatz matrices i have found the following properties: \\\\
\textit{Each column of a standard Collatz matrix has exactly one knot number.\\\\
The first row of each standard Collatz matrix is never unbranched.\\\\
The Value n is always smaller than m. There are no m x m matrices. \\\\
Each standard Collatz matrix has a kind of mirror axis (the middle row of the matrix which is always unbranched ) to which I find the following 5 main symmetries:\\\\
SM: Singlematrix\\
UM: Uppermatrix \\
MM: Mirrormatrix\\
IMM: Inverted mirrormatrix\\
USM: Unsymmetrical matrix}\\\\
The standard Collatzmatrix 11 x 10 as example for a singlematrix. Singlematrix means that there is only one unbranched row in the  middle:\\\\
$$
\begin{pmatrix}
2&4&8&16&32&64&128&256&512&\textcolor{red}{1024}\\
6&\textcolor{red}{12}&24&48&96&192&384&768&1536&3072\\
10&20&40&80&160&\textcolor{red}{320}&640&1280&2560&5120\\
14&28&\textcolor{red}{56}&112&224&448&896&1792&3584&7168\\
18&36&72&\textcolor{red}{144}&288&576&1152&2304&4608&9216\\
\textcolor{blue}{22}&\textcolor{blue}{44}&\textcolor{blue}{88}&\textcolor{blue}{176}&\textcolor{blue}{352}&\textcolor{blue}{704}&\textcolor{blue}{1408}&\textcolor{blue}{2816}&\textcolor{blue}{5632}&\textcolor{blue}{11264}\\
26&52&104&208&416&832&1664&3328&\textcolor{red}{6656}&13312\\
30&60&120&240&480&960&1920&\textcolor{red}{3840}&7680&15360\\
\textcolor{red}{34}&68&136&272&544&1088&2176&4352&8704&17408\\
38&76&152&304&608&1216&\textcolor{red}{2432}&4864&9728&19456\\
42&84&168&336&\textcolor{red}{672}&1344&2688&5376&10752&21504\
\end{pmatrix}
$$
\newpage
The standard Collatzmatrix 31 x 5 as Example for an uppermatrix. Uppermatrix means that every row above the middle axis is unbranched with the exception of the first row:
\\\\
$$
\begin{pmatrix}
2&4&8&16&\textcolor{red}{32}\\
\textcolor{blue}{6}&\textcolor{blue}{12}&\textcolor{blue}{24}&\textcolor{blue}{48}&\textcolor{blue}{96}\\
\textcolor{blue}{10}&\textcolor{blue}{20}&\textcolor{blue}{40}&\textcolor{blue}{80}&\textcolor{blue}{160}\\
\textcolor{blue}{14}&\textcolor{blue}{28}&\textcolor{blue}{56}&\textcolor{blue}{112}&\textcolor{blue}{224}\\
\textcolor{blue}{18}&\textcolor{blue}{36}&\textcolor{blue}{72}&\textcolor{blue}{144}&\textcolor{blue}{288}\\
\textcolor{blue}{22}&\textcolor{blue}{44}&\textcolor{blue}{88}&\textcolor{blue}{176}&\textcolor{blue}{352}\\
\textcolor{blue}{26}&\textcolor{blue}{52}&\textcolor{blue}{104}&\textcolor{blue}{208}&\textcolor{blue}{416}\\
\textcolor{blue}{30}&\textcolor{blue}{60}&\textcolor{blue}{120}&\textcolor{blue}{240}&\textcolor{blue}{480}\\
\textcolor{blue}{34}&\textcolor{blue}{68}&\textcolor{blue}{136}&\textcolor{blue}{272}&\textcolor{blue}{544}\\
\textcolor{blue}{38}&\textcolor{blue}{76}&\textcolor{blue}{152}&\textcolor{blue}{304}&\textcolor{blue}{608}\\
\textcolor{blue}{42}&\textcolor{blue}{84}&\textcolor{blue}{168}&\textcolor{blue}{336}&\textcolor{blue}{672}\\
\textcolor{blue}{46}&\textcolor{blue}{92}&\textcolor{blue}{184}&\textcolor{blue}{368}&\textcolor{blue}{736}\\
\textcolor{blue}{50}&\textcolor{blue}{100}&\textcolor{blue}{200}&\textcolor{blue}{400}&\textcolor{blue}{800}\\
\textcolor{blue}{54}&\textcolor{blue}{108}&\textcolor{blue}{216}&\textcolor{blue}{432}&\textcolor{blue}{864}\\
\textcolor{blue}{58}&\textcolor{blue}{116}&\textcolor{blue}{232}&\textcolor{blue}{464}&\textcolor{blue}{928}\\
\textcolor{blue}{62}&\textcolor{blue}{124}&\textcolor{blue}{248}&\textcolor{blue}{496}&\textcolor{blue}{992}\\
66&132&264&\textcolor{red}{528}&1056\\
70&140&\textcolor{red}{280}&560&1120\\
\textcolor{blue}{74}&\textcolor{blue}{145}&\textcolor{blue}{296}&\textcolor{blue}{592}&\textcolor{blue}{1184}\\
78&\textcolor{red}{156}&312&624&1248\\
\textcolor{blue}{82}&\textcolor{blue}{164}&\textcolor{blue}{328}&\textcolor{blue}{656}&\textcolor{blue}{1312}\\
\textcolor{blue}{86}&\textcolor{blue}{172}&\textcolor{blue}{344}&\textcolor{blue}{688}&\textcolor{blue}{1376}\\
\textcolor{blue}{90}&\textcolor{blue}{180}&\textcolor{blue}{360}&\textcolor{blue}{720}&\textcolor{blue}{1440}\\
\textcolor{red}{94}&188&376&752&1504\\
\textcolor{blue}{98}&\textcolor{blue}{196}&\textcolor{blue}{392}&\textcolor{blue}{784}&\textcolor{blue}{1568}\\
\textcolor{blue}{102}&\textcolor{blue}{204}&\textcolor{blue}{408}&\textcolor{blue}{816}&\textcolor{blue}{1632}\\
\textcolor{blue}{106}&\textcolor{blue}{212}&\textcolor{blue}{424}&\textcolor{blue}{848}&\textcolor{blue}{1696}\\
\textcolor{blue}{110}&\textcolor{blue}{220}&\textcolor{blue}{440}&\textcolor{blue}{880}&\textcolor{blue}{1760}\\
\textcolor{blue}{114}&\textcolor{blue}{228}&\textcolor{blue}{456}&\textcolor{blue}{912}&\textcolor{blue}{1824}\\
\textcolor{blue}{118}&\textcolor{blue}{236}&\textcolor{blue}{472}&\textcolor{blue}{944}&\textcolor{blue}{1888}\\
\textcolor{blue}{122}&\textcolor{blue}{244}&\textcolor{blue}{488}&\textcolor{blue}{976}&\textcolor{blue}{1952}\\
\end{pmatrix}
$$
\newpage
The standard Collatzmatrix 17 x 8 as a mirrormatrix. It means that the property (branched or unbranched) of each row is mirrored to the other side of the axis. For example: If the second row below the axis is unbranched then the second row above the axis is also unbranched:
\\\\
$$
\begin{pmatrix}
2&4&8&16&32&64&128&\textcolor{red}{256}\\
\textcolor{blue}{6}&\textcolor{blue}{12}&\textcolor{blue}{24}&\textcolor{blue}{48}&\textcolor{blue}{96}&\textcolor{blue}{192}&\textcolor{blue}{384}&\textcolor{blue}{768}\\
\textcolor{blue}{10}&\textcolor{blue}{20}&\textcolor{blue}{40}&\textcolor{blue}{80}&\textcolor{blue}{160}&\textcolor{blue}{320}&\textcolor{blue}{640}&\textcolor{blue}{1280}\\
\textcolor{blue}{14}&\textcolor{blue}{28}&\textcolor{blue}{56}&\textcolor{blue}{112}&\textcolor{blue}{224}&\textcolor{blue}{448}&\textcolor{blue}{896}&\textcolor{blue}{1792}\\
\textcolor{red}{18}&36&72&144&288&576&1152&2304\\
\textcolor{blue}{22}&\textcolor{blue}{44}&\textcolor{blue}{88}&\textcolor{blue}{176}&\textcolor{blue}{352}&\textcolor{blue}{704}&\textcolor{blue}{1408}&\textcolor{blue}{2816}\\
26&\textcolor{red}{52}&104&208&416&832&1664&3328\\
30&60&\textcolor{red}{120}&240&480&960&1920&3840\\
\textcolor{blue}{34}&\textcolor{blue}{68}&\textcolor{blue}{136}&\textcolor{blue}{272}&\textcolor{blue}{544}&\textcolor{blue}{1088}&\textcolor{blue}{2176}&\textcolor{blue}{4352}\\
38&76&152&304&608&1216&\textcolor{red}{2432}&4864\\
42&84&168&336&672&\textcolor{red}{1344}&2688&5376\\
\textcolor{blue}{46}&\textcolor{blue}{92}&\textcolor{blue}{184}&\textcolor{blue}{368}&\textcolor{blue}{736}&\textcolor{blue}{1472}&\textcolor{blue}{2944}&\textcolor{blue}{5888}\\
50&100&200&400&\textcolor{red}{800}&1600&3200&6400\\
\textcolor{blue}{54}&\textcolor{blue}{108}&\textcolor{blue}{216}&\textcolor{blue}{432}&\textcolor{blue}{864}&\textcolor{blue}{1728}&\textcolor{blue}{3456}&\textcolor{blue}{6912}\\
\textcolor{blue}{58}&\textcolor{blue}{116}&\textcolor{blue}{232}&\textcolor{blue}{464}&\textcolor{blue}{928}&\textcolor{blue}{1856}&\textcolor{blue}{3712}&\textcolor{blue}{7424}\\
\textcolor{blue}{62}&\textcolor{blue}{124}&\textcolor{blue}{248}&\textcolor{blue}{496}&\textcolor{blue}{992}&\textcolor{blue}{1984}&\textcolor{blue}{3968}&\textcolor{blue}{7936}\\
66&132&264&\textcolor{red}{528}&1056&2112&4224&8448\\
\end{pmatrix}
$$
\newpage
The standard Collatzmatrix 23 x 11 as an inverted mirrormatrix. The contrary of a mirrormatrix is the inverted mirrormatrix. If the second row below the axis is unbranched then the second row above the matrix is not unbranched:
\\\\
$$
\begin{pmatrix}
2&4&8&16&32&64&128&256&512&1024&\textcolor{red}{2048}\\
6&12&\textcolor{red}{24}&48&96&192&384&768&1536&3072&6144\\
\textcolor{blue}{10}&\textcolor{blue}{20}&\textcolor{blue}{40}&\textcolor{blue}{80}&\textcolor{blue}{160}&\textcolor{blue}{320}&\textcolor{blue}{640}&\textcolor{blue}{1280}&\textcolor{blue}{2560}&\textcolor{blue}{5120}&\textcolor{blue}{10240}\\
\textcolor{blue}{14}&\textcolor{blue}{28}&\textcolor{blue}{56}&\textcolor{blue}{112}&\textcolor{blue}{224}&\textcolor{blue}{448}&\textcolor{blue}{896}&\textcolor{blue}{1792}&\textcolor{blue}{3584}&\textcolor{blue}{7168}&\textcolor{blue}{14336}\\
18&36&72&144&288&\textcolor{red}{576}&1152&2304&4608&9216&18432\\
\textcolor{blue}{22}&\textcolor{blue}{44}&\textcolor{blue}{88}&\textcolor{blue}{176}&\textcolor{blue}{352}&\textcolor{blue}{704}&\textcolor{blue}{1408}&\textcolor{blue}{2816}&\textcolor{blue}{5632}&\textcolor{blue}{11264}&\textcolor{blue}{22528}\\
26&52&104&\textcolor{red}{208}&416&832&1664&3328&6656&13312&26624\\
\textcolor{blue}{30}&\textcolor{blue}{60}&\textcolor{blue}{120}&\textcolor{blue}{240}&\textcolor{blue}{480}&\textcolor{blue}{960}&\textcolor{blue}{1920}&\textcolor{blue}{3840}&\textcolor{blue}{7680}&\textcolor{blue}{15360}&\textcolor{blue}{30720}\\
\textcolor{blue}{34}&\textcolor{blue}{68}&\textcolor{blue}{136}&\textcolor{blue}{272}&\textcolor{blue}{544}&\textcolor{blue}{1088}&\textcolor{blue}{2176}&\textcolor{blue}{4352}&\textcolor{blue}{8704}&\textcolor{blue}{17408}&\textcolor{blue}{34816}\\
\textcolor{blue}{38}&\textcolor{blue}{76}&\textcolor{blue}{152}&\textcolor{blue}{304}&\textcolor{blue}{608}&\textcolor{blue}{1216}&\textcolor{blue}{2432}&\textcolor{blue}{4864}&\textcolor{blue}{9728}&\textcolor{blue}{19456}&\textcolor{blue}{38912}\\
\textcolor{blue}{42}&\textcolor{blue}{84}&\textcolor{blue}{168}&\textcolor{blue}{336}&\textcolor{blue}{672}&\textcolor{blue}{1344}&\textcolor{blue}{2688}&\textcolor{blue}{5376}&\textcolor{blue}{10752}&\textcolor{blue}{21504}&\textcolor{blue}{43008}\\
\textcolor{blue}{46}&\textcolor{blue}{92}&\textcolor{blue}{184}&\textcolor{blue}{368}&\textcolor{blue}{736}&\textcolor{blue}{1472}&\textcolor{blue}{2944}&\textcolor{blue}{5888}&\textcolor{blue}{11776}&\textcolor{blue}{23552}&\textcolor{blue}{47104}\\
50&100&200&400&800&1600&3200&6400&12800&\textcolor{red}{25600}&51200\\
54&108&216&432&864&1728&3456&6912&\textcolor{red}{13824}&27648&55296\\
58&\textcolor{red}{116}&232&464&928&1856&3712&7424&14848&29696&59392\\
62&124&248&496&992&1984&3968&\textcolor{red}{7936}&15872&31744&63488\\
\textcolor{blue}{66}&\textcolor{blue}{132}&\textcolor{blue}{264}&\textcolor{blue}{528}&\textcolor{blue}{1056}&\textcolor{blue}{2112}&\textcolor{blue}{4224}&\textcolor{blue}{8448}&\textcolor{blue}{16896}&\textcolor{blue}{33792}&\textcolor{blue}{76584}\\
\textcolor{red}{70}&140&280&560&1120&2240&4480&8960&17920&35840&71680\\
\textcolor{blue}{74}&\textcolor{blue}{145}&\textcolor{blue}{296}&\textcolor{blue}{592}&\textcolor{blue}{1184}&\textcolor{blue}{2368}&\textcolor{blue}{4736}&\textcolor{blue}{9472}&\textcolor{blue}{18944}&\textcolor{blue}{37888}&\textcolor{blue}{75776}\\
78&156&312&624&1248&2496&\textcolor{red}{4992}&9984&19968&39936&79872\\
82&164&328&656&\textcolor{red}{1312}&2624&5248&10496&20992&41984&83968\\
\textcolor{blue}{86}&\textcolor{blue}{172}&\textcolor{blue}{344}&\textcolor{blue}{688}&\textcolor{blue}{1376}&\textcolor{blue}{2752}&\textcolor{blue}{5504}&\textcolor{blue}{11008}&\textcolor{blue}{22016}&\textcolor{blue}{44032}&\textcolor{blue}{88064}\\
\textcolor{blue}{90}&\textcolor{blue}{180}&\textcolor{blue}{360}&\textcolor{blue}{720}&\textcolor{blue}{1440}&\textcolor{blue}{2880}&\textcolor{blue}{5760}&\textcolor{blue}{11520}&\textcolor{blue}{23040}&\textcolor{blue}{46080}&\textcolor{blue}{ 92160}\\
\end{pmatrix}
$$
\newpage
And finally the standard Collatzmatrix 21 x 6 as an seemingly unsymmetrical matrix:
\\\\
$$
\begin{pmatrix}
2&4&8&16&32&\textcolor{red}{64}\\
\textcolor{blue}{6}&\textcolor{blue}{12}&\textcolor{blue}{24}&\textcolor{blue}{48}&\textcolor{blue}{96}&\textcolor{blue}{192}\\
\textcolor{blue}{10}&\textcolor{blue}{20}&\textcolor{blue}{40}&\textcolor{blue}{80}&\textcolor{blue}{160}&\textcolor{blue}{320}\\
\textcolor{blue}{14}&\textcolor{blue}{28}&\textcolor{blue}{56}&\textcolor{blue}{112}&\textcolor{blue}{224}&\textcolor{blue}{448}\\
\textcolor{blue}{18}&\textcolor{blue}{36}&\textcolor{blue}{72}&\textcolor{blue}{144}&\textcolor{blue}{288}&\textcolor{blue}{576}\\
\textcolor{red}{22}&44&88&176&352&704\\
\textcolor{blue}{26}&\textcolor{blue}{52}&\textcolor{blue}{104}&\textcolor{blue}{208}&\textcolor{blue}{416}&\textcolor{blue}{832}\\
\textcolor{blue}{30}&\textcolor{blue}{60}&\textcolor{blue}{120}&\textcolor{blue}{240}&\textcolor{blue}{480}&\textcolor{blue}{960}\\
\textcolor{blue}{34}&\textcolor{blue}{68}&\textcolor{blue}{136}&\textcolor{blue}{272}&\textcolor{blue}{544}&\textcolor{blue}{1088}\\
\textcolor{blue}{38}&\textcolor{blue}{76}&\textcolor{blue}{152}&\textcolor{blue}{304}&\textcolor{blue}{608}&\textcolor{blue}{1216}\\
\textcolor{blue}{42}&\textcolor{blue}{84}&\textcolor{blue}{168}&\textcolor{blue}{336}&\textcolor{blue}{672}&\textcolor{blue}{1344}\\
46&92&184&368&\textcolor{red}{736}&1472\\
50&100&200&\textcolor{red}{400}&800&1600\\
\textcolor{blue}{54}&\textcolor{blue}{108}&\textcolor{blue}{216}&\textcolor{blue}{432}&\textcolor{blue}{864}&\textcolor{blue}{1728}\\
58&116&\textcolor{red}{232}&464&928&1856\\
\textcolor{blue}{62}&\textcolor{blue}{124}&\textcolor{blue}{248}&\textcolor{blue}{496}&\textcolor{blue}{992}&\textcolor{blue}{1984}\\
\textcolor{blue}{66}&\textcolor{blue}{132}&\textcolor{blue}{264}&\textcolor{blue}{528}&\textcolor{blue}{1056}&\textcolor{blue}{2112}\\
\textcolor{blue}{70}&\textcolor{blue}{140}&\textcolor{blue}{280}&\textcolor{blue}{560}&\textcolor{blue}{1120}&\textcolor{blue}{2240}\\
74&\textcolor{red}{148}&296&592&1184&2386\\
\textcolor{blue}{78}&\textcolor{blue}{156}&\textcolor{blue}{312}&\textcolor{blue}{624}&\textcolor{blue}{1248}&\textcolor{blue}{2496}\\
\textcolor{blue}{82}&\textcolor{blue}{164}&\textcolor{blue}{328}&\textcolor{blue}{656}&\textcolor{blue}{1312}&\textcolor{blue}{2624}\\
\end{pmatrix}
$$
\\\\
Table 1 is a summary of the first symmetries and matrices. It can be seen that there are integer sequences (A003602, A001511, A005408, A002326). \cite{4}
\newpage
\begin{table}
\caption{Summary Symmetries}
\begin{scriptsize}
\begin{center}
\begin{tabular}{|c|c|c|c|c|c|c|c|c|c|c|c|c|}
\hline
\textbf{Algo.} & \multicolumn{3}{|c|}{\textbf{Standard Collatz ma.}}& \ & \multicolumn{3}{|c|}{\textbf{Little Collatz ma.}}& \ & \multicolumn{3}{|c|}{\textbf{Big Collatz ma.}}&\textbf{Sym.}\\ \hline
\ & A5408& \ & A2326 & \ & A3602 & \ & A1511 & \ & A16754 & \ & & \ \\ \hline
$f_a(n)$ & $m_C$ &x& $n_C$ & \ & $m_L$ &x& $n_L$ & \ &$ m_B$ &x& $n_B$& \ \\ \hline
$f_1(n)$&1 &x& 1 & \ & 1 &x& 1 & \ & 1 &x& 1&UM\\ \hline
$f_3(n)$&3 &x& 2 & \ & 1 &x& 2 & \ & 9 &x& 6&UM\\ \hline
$f_5(n)$&5 &x& 4 & \ & 2 &x& 1 &  \ &25 &x& 20&SM\\ \hline
$f_7(n)$&7 &x& 3 & \ & 1 &x& 3 &  \ &49 &x& 21&UM\\ \hline
$f_9(n)$&9 &x& 6 & \ & 3 &x& 1 & \ & 81 &x& 54&MM\\ \hline
$f_{11}(n)$&11 &x& 10 & \ & 2 &x& 2 & \ & 121 &x& 110&SM\\ \hline
$f_{13}(n)$&13 &x& 12 & \ & 4 &x& 1 & \ & 169 &x& 156&SM\\ \hline
$f_{15}(n)$&15 &x& 4 &  \ &1 &x& 4 & \ & 225 &x& 60&UM\\ \hline
$f_{17}(n)$&17 &x& 8 & \ & 5 &x& 1 & \ & 289 &x& 136&MM\\ \hline
$f_{19}(n)$&19 &x& 18 & \ & 3 &x& 2 & \ & 361 &x& 342&SM\\ \hline
$f_{21}(n)$&21 &x& 6 & \ & 6 &x& 1 & \ & 441 &x& 126&USM\\ \hline
$f_{23}(n)$&23 &x& 11 & \ & 2 &x& 3 & \ & 529 &x& 253&IMM\\ \hline
$f_{25}(n)$&25 &x& 20 & \ & 7 &x& 1 & \ & 625 &x& 500&MM\\ \hline
$f_{27}(n)$&27 &x& 18 & \ & 4 &x& 2 & \ & 729 &x& 486&MM\\ \hline
$f_{29}(n)$&29 &x& 28 &  \ &8 &x& 1 & \ & 841 &x& 812&SM\\ \hline
$f_{31}(n)$&31 &x& 5 & \ & 1 &x& 5 & \ & 961 &x& 155&UM\\ \hline
$f_{33}(n)$&33 &x& 10 & \ & 9 &x& 1 & \ & 1089 &x& 330&MM\\ \hline
$f_{35}(n)$&35 &x& 12 & \ & 5 &x& 2 &  \ &1225 &x& 420&USM\\ \hline
$f_{37}(n)$&37 &x& 36 & \ & 10 &x& 1 & \ & 1369 &x& 1332&SM\\ \hline
$f_{39}(n)$&39 &x& 12 & \ & 3 &x& 3 & \ & 1521 &x& 468&USM\\ \hline
$f_{41}(n)$&41 &x& 20 & \ & 11 &x& 1 & \ & 1681 &x& 820&MM\\ \hline
$f_{43}(n)$&43 &x& 14 & \ & 6 &x& 2 & \ & 1849 &x& 602&MM\\ \hline
$f_{45}(n)$&45 &x& 12 & \ & 12 &x& 1 & \ & 2025 &x& 540&USM\\ \hline
$f_{47}(n)$&47 &x& 23 & \ & 2 &x& 4 & \ & 2209 &x& 1081&IMM\\ \hline
$f_{49}(n)$&49 &x& 21 & \ & 13 &x& 1 & \ & 2401 &x& 1029&USM\\ \hline
$f_{51}(n)$&51 &x& 8 & \ & 7 &x& 2 & \ & 2601 &x& 408&USM\\ \hline
$f_{53}(n)$&53 &x& 52 & \ & 14 &x& 1 & \ & 2809 &x& 2756&SM\\ \hline
$f_{55}(n)$&55 &x& 20 & \ & 4 &x& 3 & \ & 3025 &x& 1100&USM\\ \hline
$f_{57}(n)$&57 &x& 18 & \ & 15 &x& 1 & \ & 3249 &x& 1026&MM\\ \hline
$f_{59}(n)$&59 &x& 58 & \ & 8 &x& 2 & \ & 3481 &x& 3422&SM\\ \hline
$f_{61}(n)$&61 &x& 60 &  \ &16 &x& 1 & \ & 3721 &x& 3660&SM\\ \hline
\end{tabular}
\end{center}
\end{scriptsize}
\end{table}
\clearpage
	\subsection{Symmetry and prime numbers}
There are some interesting correlations between symmetry and prime numbers: \cite{4}\\\\ 
Each singlematrix is seemingly a prime. ( A001122 )\\\\
Each inverted mirrormatrix is seemingly a prime.  ( A139035 )\\\\
Each uppermatrix is seemingly a Mersenne number. ( A000225 ).\\\\
\section{The Collatz matrix conjecture}
	\subsection{The Collatz matrix conjecture}
The Collatz matrix conjecture establishes a connection between the standard Collatz matrix and primes like this:\\\\
\textbf{Conjecture 1. The Collatz matrix conjecture}
$$(m_C-1) \ mod \ n_C=0, \ m_C \ is \ a \ prime \ or \ a \ fermat \ pseudoprime \ base \ 2$$
$$(m_C-1) \ mod \ n_C\neq0, \ m_C \ is \ not \ a \ prime $$\\\\
$m_C=m$-value of the standard Collatz matrix\\\\
$n_C=n$-value of the standard Collatz matrix
\newpage
	\subsection{The rank of primes}
\textbf{Definition 4.1 The rank of primes}\\
With the Collatz matrix conjecture it is possible to classify primes by a rank: $$p_r=\frac{m_C-1}{n_C}$$\\\\
$p_r=$ rank of a prime\\\\
$m_C=m$-value of the standard Collatz matrix\\\\
$n_C=n$-value of the standard Collatz matrix\\\\
As can be seen in Table 2 all primes are classified by rank in the right column. In Table 3 the frequency of this rank 1-18 is listed.\\\\\begin{table}
\caption{Rank of a prime}
\begin{center}
\begin{tabular}{|c|c|c|c|c|}
\hline
\ & \multicolumn{3}{|c|}{\textbf{Collatz matrix}}& \textbf{Rank}\\ \hline
\ & A5408 & \ & A002326 & \ \\ \hline
\ & $m_C$ &x& $n_C$ & $p_r$ \\ \hline
\ &1 &x& 1 & 0\\ \hline
Prime &3 &x& 2 & 1\\ \hline
Prime &5 &x& 4 & 1\\ \hline
Prime &7 &x& 3 & 2\\ \hline
\ &9 &x& 6 & 1,333\\ \hline
Prime &11 &x& 10 & 1\\ \hline
Prime &13 &x& 12 & 1\\ \hline
\ &15 &x& 4 &  3,5\\ \hline
Prime &17 &x& 8 & 2\\ \hline
Prime &19 &x& 18 & 1\\ \hline
\ &21 &x& 6 & 3,33\\ \hline
Prime &23 &x& 11 & 2\\ \hline
\ &25 &x& 20 & 1,2\\ \hline
\ &27 &x& 18 & 1,444\\ \hline
Prime &29 &x& 28 &  1\\ \hline
Prime &31 &x& 5 & 6\\ \hline
\ &33 &x& 10 & 3,2\\ \hline
\ &35 &x& 12 & 2,83\\ \hline
Prime &37 &x& 36 & 1\\ \hline
\ &39 &x& 12 & 3,16\\ \hline
Prime &41 &x& 20 & 2\\ \hline
Prime &43 &x& 14 & 3\\ \hline
\ &45 &x& 12 & 3,66\\ \hline
Prime &47 &x& 23 & 2\\ \hline
\ &49 &x& 21 & 2,285\\ \hline
\ &51 &x& 8 & 6,25\\ \hline
\end{tabular}
\end{center}
\end{table}
\newpage
\begin{table}
\caption{Primefrequency by rank (1-18) up to number 1,000,000}
\begin{center}
\begin{tabular}{|c|c|c|c|c|}
\hline
Rank & Frequency& \ &Rank & Frequency \\ \hline
1 & 29341& \ & 10 & 1089\\ \hline
2 & 22092& \ & 11 & 278\\ \hline
3 & 5233& \ & 12 & 628\\ \hline
4 & 3655& \ & 13 & 195\\ \hline
5 & 1477& \ & 14 & 547\\ \hline
6 & 3931& \ & 15 & 248\\ \hline
7 & 694& \ & 16 & 686\\ \hline
8 & 2781& \ & 17 & 115\\ \hline
9 & 579& \ & 18 & 432\\ \hline
\end{tabular}
\end{center}
\end{table}
	\subsection{Rank and symmetry}
With the rank of a prime it is possible to make propositions about the symmetry of the Collatz matrices:\\\\
\textit{Singlematrix:\\\\
Each Collatz matrix with the rank 1 are singlematrices.\\\\
Inverted mirrormatrix:\\\\
Each Collatz matrix with the rank 2 and an odd n-value of the Collatzmatrix are inverted mirrormatrices.\\\\
Uppermatrices:\\\\
Each Collatz matrix with 1 as an m-value of the little Collatzmatrix are uppermatrices.}
\clearpage
	\subsection{Fermat pseudoprimes (base 2) classified by the little Collatz matrix}	
As mentioned in 4.1, the Collatz matrix conjecture shows a necessary criterion for a prime. The numbers which have this criterion but are no primes are Fermat pseudoprimes to base 2. With the little Collatz matrix it is now possible to get a new classification for this pseudoprimes (Table 4). As can be seen there are many pseudos with 1, a few with 2 and so on... 
\\\\
\begin{table}
\caption{Classified Pseudoprimes}
\begin{scriptsize}
\begin{center}
\begin{tabular}{|c|c|c|c|c|c|c|c|c|}
\hline
\multicolumn{3}{|c|}{\textbf{Collatz matrix}}& \ & \multicolumn{3}{|c|}{\textbf{Little Collatz matrix}}&\textbf{Rank}\\ \hline
341 &x& 10 & \ & 86 &x& 1 & 34\\ \hline
561 &x& 40 & \ & 141 &x& 1 &14\\ \hline
645 &x& 28 & \ & 162 &x& 1 &23\\ \hline
1105 &x& 24 & \ & 277 &x& 1&46\\ \hline
1387 &x& 18 & \ & 174 &x& 2&77\\ \hline
1729 &x& 36 & \ & 433 &x& 1&48\\ \hline
1905 &x& 28 & \ & 477 &x& 1&68\\ \hline
2047 &x& 11 &  \ &1 &x& 11&186\\ \hline
2465 &x& 56 & \ & 617 &x& 1&44\\ \hline
2701 &x& 36 & \ & 676 &x& 1&75\\ \hline
2821 &x& 60 & \ & 706 &x& 1&47\\ \hline
3277 &x& 28 & \ & 820 &x& 1&117\\ \hline
4033 &x& 36 & \ & 1009 &x& 1&112\\ \hline
4369 &x& 16 & \ & 1093 &x& 1&273\\ \hline
4371 &x& 230 &  \ &547 &x& 2&19\\ \hline
4681 &x& 15 & \ & 1171 &x& 1&312\\ \hline
5461 &x& 14 & \ & 1366 &x& 1&390\\ \hline
\end{tabular}
\end{center}
\end{scriptsize}
\end{table}
\clearpage
\section{Kaiser's conjecture}
	\subsection{Mersenne primes}
A Mersenne number is a number of the form $2^n-1$. Mersenne primes are Mersenne numbers which are also primes. The exponents n of such a Mersenne prime are 2,3,5,7,13,17,19,31,61... (A000043) \cite{4,5}
	\subsection{Kaiser’s conjecture}
When comparing Collatz matrices with Mersenne prime exponents the following context can be seen:\\\\
\textbf{Conjecture 2. Kaiser's conjecture}\\
The exponent n of a Mersenne prime of the form $2^n-1$ has to be singular. Singular means that it exists only one Collatz matrix with this exponent n as n-value $n_C$.\\\\
For example:\\\\
There is one Collatz matrix with the value $n_C=3$. It is 7x3...\\\\
There is one Collatz matrix with the value $n_C=5$. It is 31x5...\\\\
There is one Collatz matrix with the value $n_C=7$. It is 127x7...\\\\
But there are 3 Collatzmatrices with the value $n_C=11$! These Collatz matrices are 23x11,89x11 and 2047x11. So 11 is not a Mersenne prime exponent.
\newpage 
In Table 5 you can see the frequency of the n-values $n_C$ 1-19 for Collatz matrices up to the Collatz matrix 1.999.999 x 6440.\\\\ 
\begin{table}
\caption{Frequency n-value}
\begin{center}
\begin{tabular}{|c|c|}
\hline
$n_C$ & Frequency\\ \hline
1 & 1\\ \hline
2 & 1\\ \hline
3 & 1\\ \hline
4 & 2\\ \hline
5 & 1\\ \hline
6 & 3\\ \hline
7 & 1\\ \hline
8 & 4\\ \hline
9 & 2\\ \hline
10 & 5\\ \hline
11 & 3\\ \hline
12 & 16\\ \hline
13 & 1\\ \hline
14 & 5\\ \hline
15 & 5\\ \hline
16 & 8\\ \hline
17 & 1\\ \hline
18 & 24\\ \hline
19 & 1\\ \hline
\end{tabular}
\end{center}
\end{table}\\\\
\section{Summary}
I have come to the conclusion that there is a deep relationship between the Collatz matrices and primes . As can be seen there is a way to test Mersenne primes theoreticaly. Furthermore, there are many structures and symmetries to find in these Collatz matrices. \\\\
\section{Acknowledgement}
Thanks to the internet and to the books of Lothar Papula.

\end{document}